\documentclass[12pt]{article}
\usepackage{amsmath,eufrak}
\usepackage{amsthm}

\newcommand{\aaa}{\mathcal{A}} %basic_algebra
\newcommand{\be}[1]{\begin{equation} \label{#1} }
\newcommand{\bra}[1]{\pmb{\langle}#1\pmb{|}}
\newcommand{\braket}[2]{\pmb{\langle}#1\pmb{\mid}#2\pmb{\rangle}} 
\newcommand{\borderop}{\delta} 
\DeclareMathOperator{\cardinality}{card}

\newcommand{\cdt}{,\ldots,}
\newcommand{\dff}{\itshape\sffamily} %font_for_definitions
\newcommand{\dm}{\mathcal{D}} %diff_module
\newcommand{\DMCA}{\mathfrak{DM}}
\newcommand{\hhh}{\mathcal{H}} %Hilbert_space

\newcommand{\ia}{\Omega} %incidence_algebra
\newcommand{\iii}{\mathfrak{I}} %diff_ideal
\newcommand{\kd}{\text{\large\texttt{d}}}%Kaehler_diff
\newcommand{\ket}[1]{\pmb{|}#1\pmb{\rangle}}%ket-vector

\newcommand{\kkk}{\mathcal{K}}%posets_complexes
\DeclareMathOperator*{\lspan}{span}%linear_span
\newcommand{\nbd}[1]{\(#1\)\nobreakdash-\hspace{0pt}}
\newcommand{\phio}{\ensuremath{\overline\phi}}
\newcommand{\pp}{\ensuremath{\pi}}%simplicial_mapping

\newcommand{\restr}[2]{{\bigl.#1\bigr\vert}_{#2}}%restrict_mapping
\newcommand{\rf}[1]{(\ref{#1})}%reference_in_()
\newcommand{\SC}{\ensuremath{\mathfrak{SC}}}%category_of_complexes

\newcommand{\story}[1]{\mbox{\ensuremath{\pmb\prec\hspace{-.3em}
#1\hspace{-.3em}\pmb\succ}}} 

\newcommand{\vvv}{\ensuremath{V}}%set_of_vertices

\newtheorem{dfn}{Definition}
\newtheorem{theorem}{Theorem}
\newtheorem{lemma}[theorem]{Lemma}

\newcommand{\bracket}[3]{\bra{#1}#2\ket{#3}} 
\newcommand{\iak}{\ensuremath{\ia(\kkk)}} 
\newcommand{\ketbra}[2]{\ket{#1}\bra{#2}} 
\newcommand{\pkk}{\ensuremath{\pp:\kkk'\to\kkk}}
\newcommand{\ude}{\ensuremath{\ia_u}} %univ_diff_envelope
\newcommand{\ukd}{\ensuremath{\kd_u}} %univ_Kaehler

\author{Roman R. Zapatrin}
\title{Incidence algebras of simplicial complexes}

\begin{document}

\maketitle

\begin{center}
Department of Mathematics, SPb UEF, Griboyedova 30/32, \\
191023, St-Petersburg, Russia
\end{center}

\begin{abstract}
With any locally finite partially ordered set $\kkk$ its incidence
algebra $\iak$ is associated. We shall consider algebras over fields
with characteristic zero. In this case there is a correspondence
$\kkk \leftrightarrow \iak$ such that the poset $\kkk$ can be
reconstructed from its incidence algebra up to an isomorphism ---
due to Stanley theorem. In the meantime, a monotone mapping between
two posets in general induces no homomorphism of their incidence
algebras.

In this paper I show that if the class of posets is confined to
simplicial complexes then their incidence algebras acquire the
structure of differential moduli and the correspondence 
$\kkk\leftrightarrow\ia$ is a contravariant functor.
\end{abstract}

\section*{Introduction}\label{s1}

This paper, being motivated by physical problems, brings together 
the issues which traditionally belong to `disjoint' areas of 
mathematics:  combinatorics and differential moduli and, besides 
that, uses the notation from quantum mechanics.

It is shown that simplicial complexes resemble differential
manifolds from the algebraic point if view, namely, their incidence
algebras are similar to algebras of exterior differential forms on
manifolds: they are graded, and possess an analog of Cartan
differential.

To make the paper self-consistent, I begin it with an outline of
basic definitions and results.

\paragraph{Dirac notation.} Let $\hhh$ be a finite-dimensional linear
space with a basis labelled by an index set $\kkk$, and $\hhh^*$
be its dual. Write down the elements of $\hhh$ as

\[
h\in \hhh \;\Leftrightarrow\; h = \sum_{P\in\kkk}c_P\ket{P}
\]

Since the dimension of $\hhh$ is finite, the same index set $\kkk$ is
used to label the dual basis in the space $\hhh^*$

\[
h^*\in \hhh^*\; \Leftrightarrow\; h^* = \sum_{P\in\kkk}c_P\bra{P}
\]

\noindent such that

\be{e50}
\braket{P}{Q} = \delta_{PQ} =
\begin{cases}
1, & \text{if $P=Q$}\\
0  & \text{otherwise}
\end{cases}
\end{equation}

The elements of $\hhh^*$ are called bra-vectors and the elements of
$\hhh$ are ket-vectors (the terms derived from splitting
the word `bracket').

Let $A:\hhh\to\hhh$ be a linear operator and $A^*:\hhh^*\to\hhh^*$
be its adjoint. In the dual bases we have for any $P,Q \in\kkk$:

\be{e50a}
\bigl(A^*(\bra{P})\bigr)\ket{Q} = \bra{P}\bigl(A(\ket{Q})\bigr)
\end{equation}

\noindent therefore with no confusion we can use the same symbol,
say, $a$ for both $A$ and its adjoint $A^*$: 

\[
\begin{array}{crcl}
A:&\ket{Q}&\mapsto&a\ket{Q}\\
A^*:&\bra{P}&\mapsto&\bra{P}a
\end{array}
\]

\noindent and the identity \rf{e50a} reads

\[
(\bra{P}a)\ket{Q}=
\bra{P}(a\ket{Q})=
\bracket{P}{a}{Q}=
a_{pq}
\]

Then $A$ can be written down as

\[
A = \sum_{P,Q\in\kkk}a_{PQ}\ketbra{P}{Q}
\] 

\noindent and the product is calculated in accordance with
\rf{e50}:

\be{e51}
A\cdot B =
\left(\sum_{P,Q\in\kkk}a_{PQ}\ketbra{P}{Q}\right)\cdot
\left(\sum_{R,S\in\kkk}b_{RS}\ketbra{R}{S}\right) =
\sum_{P,Q,S\in\kkk}a_{PQ}b_{QS}\ketbra{P}{S}
\end{equation}

This notation was introduced by P.A.M.Dirac \cite{dirac} for state
vectors in quantum mechanics.

\paragraph{Incidence algebras.} The Dirac notation turns out to be
natural for incidence algebras. Let $\kkk$ be an arbitrary finite
poset. Denote by $\hhh$ its linear span

\[
\hhh = \lspan\{\ket{P}:\, P\in\kkk\} =
\left\{\sum_{P\in\kkk}c_P\ket{P}\right\}
\]

\noindent with the coefficients taken from a field with
characteristic zero.

\begin{dfn} The {\dff incidence algebra} of a poset $\kkk$ is the
following linear span

\be{e51a}
\ia = \iak =
\lspan\{ \ketbra{P}{Q}:\,P,Q\in\kkk\;\text{\textup{and}}\;P\le Q\}
\end{equation}

\noindent and the product defined on the basic elements according
to \rf{e51a}:

\[
\ketbra{P}{Q}\cdot\ketbra{R}{S} =
\ket{P}\braket{Q}{R}\bra{S} =
\delta_{QR}\ketbra{P}{S}
\]

\end{dfn}

This definition of product is correct due to the transitivity of
partial orders:

\[
\ketbra{P}{Q},\,\ketbra{Q}{S}\in\ia
\quad\Rightarrow\quad
P\le Q\;\text{and}\;Q\le S
\quad\Rightarrow\quad
P\le S
\quad\Rightarrow\quad
\ketbra{P}{S}\in\ia
\]

Incidence algebras were introduced by Rota \cite{rota}. It was
proved by Stanley \cite{stanley} that a poset $\kkk$ can be
reconstructed from its incidence algebra $\iak$ up to a poset
isomorphism.

\medskip

Meanwhile, poset homomorphisms (namely, monotone mappings) induce
no homomorphism of incidence algebras. However, and this is the
contents of this paper, the situation drastically changes when the
class of posets is restricted to simplicial complexes.

\paragraph{The category $\SC$ of simplicial complexes.} For the sake
of self-consistency I give a brief account of the standard theory,
mainly to introduce the notation. Let $\vvv$ be a non-empty finite
set, call the elements of $\vvv$ {\dff vertices}.

\begin{dfn}
A collection $\kkk$ of non-empty subsets of $\vvv$ is called (abstract)
{\dff simplicial complex} with the set of vertices $\vvv$ whenever

\begin{itemize}
\item $\forall v\in\vvv \quad \{v\}\in\kkk$
\item $\forall P\in\kkk,\, \forall Q\subseteq\vvv \quad
Q\subseteq P\Rightarrow Q\in\kkk$
\end{itemize}

\end{dfn}

Evidently, $\kkk$ is a poset with respect to set inclusion. The 
elements of $\kkk$ are called {\dff simplices}.

\begin{dfn}
Let $\kkk$, $\kkk'$ be two simplicial complexes with the sets of
vertices $\vvv$, $\vvv'$, respectively. A mapping $\pkk$ is called
{\dff simplicial} if

\begin{itemize}
\item vertices are mapped on vertices: $\quad\vvv'\pp\subseteq\vvv$
\item $\restr{\pp}{\vvv'}$ completely determines $\pp$ on the whole 
 $\kkk'$: $\quad\forall P'\in\kkk'\quad P'\pp = \cup_{v'\in 
 P'}\{v'\}\pp$ 
\end{itemize}

\end{dfn}

\noindent \textbf{Remark.} Sometimes (in particular, in the rest of 
this paper) the notion of simplicial mapping is referred to a 
mapping $\pp$ between vertices, then the second condition reads:

\be{e80}
\{v_0'\cdt v_n'\}\in\kkk' \quad\Rightarrow\quad
\{v_0'\pp\cdt v_n'\pp\}\in\kkk'
\end{equation}

\noindent For instance, let
\newcounter{NN}
\setcounter{NN}{0}
\(
\kkk' =
\setlength{\unitlength}{.1ex}
\begin{picture}(80,30)
\multiput(5,7)(30,0){3}{\circle*{6}}
\multiput(5,7)(30,0){2}{\line(1,0){26}}
\multiput(-7,15)(30,0){3}{\addtocounter{NN}{1}
\mbox{\scriptsize{\theNN '}}}
\end{picture}
\)
and
\setcounter{NN}{0}
\(
\kkk =
\setlength{\unitlength}{.1ex}
\begin{picture}(80,30)
\multiput(5,7)(30,0){3}{\circle*{6}}
\multiput(5,7)(30,0){2}{\line(1,0){26}}
\multiput(-7,15)(30,0){3}{\addtocounter{NN}{1}
\mbox{\scriptsize{\theNN}}}
\end{picture}
\)
and let $\pp_1$, $\pp_2$ be

\setlength{\unitlength}{.4ex}

\[
\begin{array}{ccc}

\begin{picture}(45,80)
    \setcounter{NN}{3}
    \multiput(0,5)(0,30){3}{\mbox{\theNN '}
    \addtocounter{NN}{-1}}
    \multiput(7,8)(0,30){3}{\circle*{3}}
    \multiput(7,10)(0,30){2}{\line(0,1){26}}

    \setcounter{NN}{3}
    \multiput(40,5)(0,30){3}{\mbox{\theNN}
    \addtocounter{NN}{-1}}
    \multiput(37,8)(0,30){3}{\circle*{3}}
    \multiput(37,10)(0,30){2}{\line(0,1){26}}

    \put(10,11){\vector(1,1){24}}
    \put(10,38){\vector(1,0){24}}
    \put(10,68){\vector(1,0){24}}
\end{picture}
&
\phantom{\text{a gap}}
&
\begin{picture}(45,80)
    \setcounter{NN}{3}
    \multiput(0,5)(0,30){3}{\mbox{\theNN '}
    \addtocounter{NN}{-1}}
    \multiput(7,8)(0,30){3}{\circle*{3}}
    \multiput(7,10)(0,30){2}{\line(0,1){26}}

    \setcounter{NN}{3}
    \multiput(40,5)(0,30){3}{\mbox{\theNN}
    \addtocounter{NN}{-1}}
    \multiput(37,8)(0,30){3}{\circle*{3}}
    \multiput(37,10)(0,30){2}{\line(0,1){26}}

    \put(10,8){\vector(1,0){24}}
    \put(10,35){\vector(1,-1){24}}
    \put(10,68){\vector(1,0){24}}
\end{picture}

\cr

\pp_1:\kkk'\to\kkk && \pp_2:\kkk'\to\kkk
\end{array}
\]

\noindent then $\pp_1$ is simplicial and $\pp_2$ is not (since
$\{1',2'\}\pp_2 = \{1,3\}\not\in\kkk$).

Denote by $\SC$ the category whose objects are simplicial complexes
and whose arrows are simplicial mappings

\[
\SC = (\text{\texttt{simplicial complexes}},\;
\text{\texttt{simplicial mappings}})
\]

\noindent It follows immediately from the definition that
simplicial mappings are monotone with respect to set inclusion and
that $\SC$ is (not a full) sub-category of the category
\(\mathfrak{POSET} = (\text{\texttt{posets}}, 
\text{\texttt{monotone\nolinebreak{ } mappings}})\)

\paragraph{Dirac notation for homological operations.} Fix an
enumeration of the vertices of $\kkk$, then for any simplex
\(P=\{v_0\cdt v_n\}\in\kkk\)
and any its vertex $v_i$ the incidence coefficient is defined

\be{e53a}
\epsilon_{v_iP} = (-1)^i
\end{equation}

\noindent A {\dff face} $P_v$ of a simplex $P$ is its subset
$P\setminus \{v\}$, we write

\be{e53m}
P_v = P-v;\; P = p_v+v
\end{equation}

\noindent The {\dff dimension} of a simplex $P$ is the number of 
its vertices minus one:

\be{e53}
\dim P = \cardinality P - 1
\end{equation}

\noindent Denote by $\kkk^n$ the \nbd{n}skeleton of $\kkk$ --- the 
set of its simplices of dimension $n$

\[
\kkk^n =
\{P\in\kkk:\:\dim P = n\}
\]

\noindent and consider the linear spans

\[
\hhh^n = \lspan\kkk^n =
\left\{\sum_{P\in\kkk^n}c_P\ket{P}\right\}
\]

\begin{dfn}
The {\dff border operator} $\borderop:\hhh^n\to\hhh^{n-1}$ acts as

\be{e53b}
\borderop\ket{P} = \sum_{v\in P}\epsilon_{vP}\ket{P_v}
\end{equation}

\noindent and then extends to the space $\hhh$ (assuming
$\borderop\ket{v}=0,\,\forall v\in\vvv$):

\[
\hhh = \oplus \hhh^n = \lspan\kkk
\]

\end{dfn}

It is proved that $\borderop^2$ is always zero. Due to Dirac
notation the same symbol $\borderop$ is used to denote its adjoint,
called {\dff coborder operator} acting from $\hhh^{*n}$ to
$\hhh^{*n+1}$, so

\be{e54e}
\forall P\in\kkk^n \qquad
\left\lbrace\begin{array}{rcl}
\bra{P}\borderop &\in& \hhh^{*n+1}\\
\borderop\ket{P} &\in& \hhh^{n-1}
\end{array}
\right.
\end{equation}

\paragraph{Finite-dimensional differential moduli.} Recall the
basic definitions concerning finite-dimensional analogs of moduli
of exterior differential forms. Let $\aaa$ be a semisimple
finite-dimensional commutative algebra.

\begin{dfn}\label{dfdm}
A {\dff differential module} $\dm$ over a basic algebra $\aaa$ is a
triple

\[
\dm = (\ia, \aaa, \kd)
\]

\noindent where $\ia$ is a graded algebra

\be{e54a}
\ia = \ia^0\oplus\ia^1\oplus\dotsb, \qquad \ia^0 = \aaa
\end{equation}

\noindent equipped with the {\dff K\"ahler differential}
$\kd:\ia^n\to\ia_{n+1}$ such that for any $\omega^r\in\ia^r$,
$\omega^s\in\ia^s$

\be{e54}
\begin{array}{l}
\kd^2 = 0 \cr
\kd(\omega^r\cdot\omega^s) =
\kd\omega^r\cdot\omega^s +
(-1)^r\omega^r\cdot\kd\omega^s
\end{array}
\end{equation}
\end{dfn}

\noindent The second equality is called {\dff graded Leibniz rule}.

\paragraph{Universal differential envelope.} Given an algebra $\aaa$,
any differential module over it can be obtained as a quotient of
a universal object $\ude = \ude(\aaa)$, called universal
differential envelope of $\aaa$ over appropriate differential ideal
$\iii$. Recall the necessary definitions (for details the Reader is
referred to \cite{kastler}).

\medskip 

Consider the tensor product $\aaa\otimes\aaa$ and define the
operator $m:\aaa\otimes\aaa\to\aaa$ as follows:

\[ m(a\otimes b) = a\cdot b \]

\noindent Then consider its kernel

\[ \ude^1 = \ker m \]

\noindent define by induction

\[ \ude^{n+1} = \ude^n\oplus_\aaa \ude^1 \]

\noindent and form the sum

\be{e55}
\ude = \bigoplus_{n=0}^\infty \ude^n
\end{equation}

Define the K\"ahler differential $\ukd$ first on $\ude^0 =
\aaa$:

\[
\ukd a = \mathbf{1} \otimes a - a\otimes\mathbf{1}
\]

\noindent and then extend it by induction to higher degrees using
the identities \rf{e54}, for instance

\be{e55a}
\ukd(a\ukd b) = \ukd a\cdot\ukd b
\end{equation}

\begin{dfn}
The {\dff universal differential envelope} $\ude = \ude(\aaa)$ of
the algebra $\aaa$ is the differential module $(\ude, \aaa, \ukd)$.
\end{dfn}

\section{Differential structure of incidence algebras}

In this section we introduce a particular representation for
universal differential envelopes of finite-dimensional algeras
called {\dff stories semantics}. Let $\kkk$ be an arbitrary (yet
structureless) finite set, denote by $\aaa$ the algebra of all
complex-valued functions on $\kkk$.

\paragraph{Stories semantics for universal differential envelope.}
Call the elements of $\kkk$ statements, and consider first all
possible sequences of statements of finite length. A {\dff 
homogeneous \nbd{n}story} is a sequence $\story{P_0\cdt P_n}$ whose 
no neighbor statements are the same. Denote by $\ia_S^n$ the linear 
span of all homogeneous \nbd{n}stories:

\be{e56}
\ia_S^n = \lspan\bigl\{\story{P_0\cdt P_n}:\;
\forall i=1\cdt n \quad
P_{i-1}\neq P_i\bigr\}
\end{equation}

\noindent and define the product of two stories as follows:

\be{e56a}
\story{P_0\cdt P_n}\cdot\story{Q_0\cdt Q_m} =
\begin{cases}
\story{P_0\cdt P_nQ_1\cdt Q_m}, &\text{if \(P_n=Q_0\)}\cr
0 &\text{otherwise}
\end{cases}
\end{equation}

\noindent which, being extended by linearity to the direct sum
\(\bigoplus_{n=0}^\infty\ia_S^n\) makes it graded algebra, call it 
{\dff stories algebra}.

\medskip 

Now let us write down the explicit form of the K\"ahler
differential. For $\story{P}\in\aaa$ we have:

\[
\ukd\story{P}=\sum_{Q\not=P}(\story{QP}-\story{PQ})
\]

\noindent and for arbitrary $\story{P_0\cdt P_n}\in\ude^n$

\begin{multline}\label{e59}
\ukd\story{P_0\cdt P_n}=
\sum_{Q:\,Q\neq P_0}\story{QP_0\cdt P_n}+\\+
\sum_{k=1}^n(-1)^k\sum_{Q:P_{k-1}\neq Q\neq P_k}
\story{P_0\cdt P_{k-1}QP_k\cdt P_n}+ \\+
(-1)^{n+1}\sum_{Q:P_n\neq Q}\story{P_0\cdt P_nQ}
\end{multline}

\begin{lemma} 
The universal differential envelope of a finite-dimensional 
semisimple algebra $\aaa$ is isomorhic (as graded algebra) to the 
stories algebra: 

\[
\ude = \bigoplus_{n=0}^\infty\ia_S^n
\]
\end{lemma}

\paragraph{Simplicial differential ideals.} Let us specify the
structure of the set $\kkk$ assuming it to be a simplicial complex.
Our goal is to find the explicit form of the differential ideal
$\iii$ giving rise to its incidence algebra $\iak$. To gather it, 
first let us classify the stories in a way taking into account that 
$\kkk$ is a simplicial complex.

\begin{dfn}
A story $\story{P_0\cdt P_n}$ is called {\dff fair} whenever
$P_{i-1} = P_i-v$ --- see \rf{e53m}. Otherwise the story
$\story{P_0\cdt P_n}$ is {\dff unfair}.
\end{dfn}

Form the linear span \(\iii_N=\oplus_{n=0}^\infty\,\iii_N^n\),
where

\be{e57}
\iii_N^n =
\lspan\{\text{\ttfamily unfair \nbd{n}stories}\} =
\lspan\{\story{P_0\cdt P_n}:\;
\lnot(\forall i\:P_{i-1}=P_i-v)\}
\end{equation}

Fix an enumeration of vertices of $\kkk$ and with each fair story $w$
associate a number $\epsilon_w=\pm1$ as follows:

\[
\epsilon_w = \prod_{i=1}^n\epsilon_{v_iP_i}
\]

\noindent where $\epsilon_{v_iP_i}$ is the incidence coefficient
\rf{e53a}. Let $P_0$, $P_n$ be two simplices such that
$P_0\subseteq P_n$, and $\dim P_0 - \dim P_n$ is exactly $n$. For
any two fair \nbd{n}stories \(w=\story{P_0P_1\cdt P_{n-1}P_n}\) and
\(w'=\story{P_0P_1'\cdt P_{n-1}'P_n}\) having common initial and
final statements form the difference \(\epsilon_ww -
\epsilon_{w'}w'\), and consider the linear hull of all such
differences

\[
\iii_S^n = \lspan_w\{
\epsilon_ww - \epsilon_{w'}w':\:
w,w'\;\text{as described above}\}
\]

\noindent Take for each $n$ the sum

\[
\iii^n = \iii_N^n\oplus\iii_S^n
\]

\noindent and form the following graded linear space

\be{e58}
\iii = \bigoplus_{n=1}^\infty\iii^n
\end{equation}

\begin{lemma}
$\iii$ is a differential ideal in $\ude$.
\end{lemma}

\begin{proof}
Evidently $\iii_N=\oplus_n\iii_N^n$ is an ideal in $\ude$. Besides
that, a product of an element of $\iii_S^n$ and a fair \nbd{m}story
is always in $\iii_S^{m+n}$, therefore $\iii$ is an ideal. It remains
to prove that the ideal $\iii$ is differential:
$\ukd(\iii^n)\in\iii^{n+1}$.

First let $w\in\iii_N^n$, consider $\ukd w$ as the sum \rf{e59}.  
The 1st and the 3rd summands of \rf{e59} are always in 
$\iii_N^{n+1}$, the same for any term from the middle sum of 
\rf{e59} with the only possible exception when \(w=\story{P_0\cdt 
P_iP_{i+1}\cdt P_n}\) such that both $\story{P_0\cdt P_i}$ and 
$\story{P_{i+1}\cdt P_n}$ are fair stories, while 
$P_i=P_{i-1}-u-v$. Then

\[
\ukd w =
\nu+(-1)^{i+1}(\rho+\rho')
\]

\noindent where $\nu$ is the sum of elements of \rf{e59} from
$\iii_N^{n+1}$, and $\rho$, $\rho'$ are the following fair stories:

\[
\begin{array}{rcl}
\rho  &=& \story{P_0\cdt P_i,P_{i+1}-u,P_{i+1}\cdt P_n}\cr
\rho' &=& \story{P_0\cdt P_i,P_{i+1}-v,P_{i+1}\cdt P_n}
\end{array}
\]

\noindent for which $\epsilon_\rho=-\epsilon_{\rho'}$, therefore
$\rho+\rho'\in\iii_S^{n+1}$ and

\[
\ukd\iii_N^n\subseteq\iii_N^{n+1}\oplus\iii_S^{n+1}
\]

Now let
\(w=\epsilon_{\rho}\rho-\epsilon_{\rho'}\rho'\in\iii_S^n\).
That is,

\[
\rho=\story{P_0,P_1\cdt P_{n-1},P_n}
,\quad
\rho'=\story{P_0,P_1'\cdt P_{n-1}',P_n}
\]

Consider the three sums \rf{e59} for \(\ukd w\). All the terms in
the 2nd sum \rf{e59} will be in $\iii_N^{n+1}$. The terms from the first
sum will also belong to $\iii_N^{n+1}$ with the only exception ---
the terms of the form
\(k=\epsilon_\tau\tau-\epsilon_{\tau'}\tau'\), where

\[
\tau=\story{P_0-v,P_0}\rho
,\quad
\tau'=\story{P_0-v,P_0}\rho'
\]

\noindent but for them
\(\epsilon_\tau=\epsilon_{vP_0}\epsilon_\rho\)
and
\(\epsilon_{\tau'}=\epsilon_{vP_0}\epsilon_\rho'\),
therefore
\(\epsilon_\tau\epsilon_{\tau'}=\epsilon_\rho\epsilon_{\rho'}\)
and $k\in\iii_S^{n+1}$, so
\(\ukd\iii_S^n\subseteq\iii_N^n\oplus\iii_S^n\). This completes the
proof:

\[
\ukd(\iii_N^n\oplus\iii_S^n)\subseteq
\iii_N^{n+1}\oplus\iii_S^{n+1}
\]

\end{proof}

\begin{theorem}
The quotient $\ude/\iii$ and the incidence algebra $\iak$ are
isomorphic algebras.
\end{theorem}

\begin{proof}
Consider the mapping $\sigma:\ude\to\ia$ defined on any story
$w=\story{P_0\cdt P_n}\in\ude^n$ as

\be{e61}
\sigma(w)=
\begin{cases}
\epsilon_w\ketbra{P_0}{P_n}, &\text{if $w$ is a fair story}\\
0 &\text{otherwise}
\end{cases}
\end{equation}

\noindent then $\ker\sigma=\iii_N\oplus\iii_S=\iii$, so
$\ia\simeq\ude/\iii$ as graded linear spaces. To verify that
$\sigma$ preserves products, let $w=\story{P_0\cdt P_m}$ and
$w'=\story{Q_0\cdt Q_n}$ be two fair stories. If $P_m\neq Q_0$
everything is trivial. Suppose $P_m=Q_0$, then
\(ww'=\story{P_0\cdt P_mQ_1\cdt Q_n}\) and

\[
\sigma(ww')=
\epsilon_{(ww')}\ketbra{P_0}{Q_n}=
\epsilon_{w}\epsilon_{w'}\ketbra{P_0}{P_m}\ketbra{Q_0}{Q_n}=
\sigma(w)\sigma(w')
\]

\end{proof}

\paragraph{The induced differential structure on $\iak$.} Having a
projection of the universal differential envelope $\ude$ onto the
incidence algebra $\iak$ we first conclude that $\iak$ necessarily has
the structure of a differential module over $\aaa$, namely, that
induced by the projection $\sigma$ onto the quotient. Now, having
the expression \rf{e61} for $\sigma$, let us explicitly calculate
the form of the differential $\kd$ on $\iak$ according to the formula

\be{eikd}
\kd=\sigma^{-1}\circ\ukd\circ\sigma
\end{equation} 

\begin{theorem}
The differential in the incidence algebra $\iak$ of a simplicial
complex $\kkk$ has the following form. Let $\ketbra{P}{Q}\in\ia^n$,
then

\be{e62}
\kd\ketbra{P}{Q}=
\ket{\borderop P}\bra{Q}-
(-1)^n\ket{P}\bra{Q\borderop}
\end{equation}

\noindent where $\borderop$ is the symbol \rf{e54e} for both border
and coborder operations.

\end{theorem}

\begin{proof}
Begin with $\ia^0=\aaa$. In this case
\(\sigma^{-1}=\restr{\text{id}}{\aaa}\), therefore

\[
\kd\ketbra{P}{P}=
\sigma(\ukd\story{P})=
\sum_{v\in P}\epsilon_{vP}\ketbra{P-v}{P}-
\sum_{u:P+u\in\kkk}\epsilon_{u(P+u)}\ketbra{P}{P+u}
\]

\noindent According to \rf{e53b}, the first term is
\(\ketbra{\borderop P}{P}\) while the second is 
\(\ketbra{P}{P\borderop}\), so

\[
\kd\,\ketbra{P}{P}=
\ket{\borderop P}\bra{P}-\ket{P}\bra{P\borderop}
\]

\noindent Now let $\ketbra{P}{Q}\in\ia^1$, then $P=Q-v$ for some
$v\in\vvv$ and

\[
\ketbra{P}{Q}=
\ketbra{P}{P}\cdot\epsilon_{vQ}\,\kd(\ketbra{Q}{Q})
\]

\noindent therefore it follows from \rf{e55a}
\(\epsilon_{vQ}=\bracket{P}{\borderop}{Q}\) that

\begin{multline*} 
\kd\ketbra{P}{Q}=
\kd\ketbra{P}{P}\cdot\epsilon_{vQ}\kd(\ketbra{Q}{Q})=
(\ket{\borderop P}\bra{P}-\ket{P}\bra{P\borderop})\epsilon_{vQ}
(\ket{\borderop Q}\bra{Q}-\ket{Q}\bra{Q\borderop})=\\
\epsilon_{vQ}(\ket{\borderop P}\bracket{P}{\borderop}{Q}\bra{Q}-
\ket{Q}\bracket{P}{\borderop}{Q}\bra{Q\borderop})=
\ket{\borderop P}\bra{Q}+\ket{P}\bra{Q\borderop}
\end{multline*} 

For higher degrees the formula \rf{e62} is proved by induction.  
First note that $\kd$ enjoys the Leibniz rule as both $\sigma$, 
$\sigma'$ in \rf{eikd} preserve products. Then represent any 
$\ketbra{P}{Q}\in\ia^n$ as a product

\[
\ketbra{P}{Q}=
\ketbra{P}{Q-v}\cdot\ketbra{Q-v}{Q}\in
\ia^{n-1}\cdot\ia^n
\]

\noindent and carry out a routine calculation.

\end{proof}

\paragraph{Summary of this section.} It was established that the
incidence algebra $\iak$ of any simplicial complex $\kkk$ is a
differential module over the algebra $\aaa$ of all functions on $\kkk$
(see Definition \ref{dfdm}). The algebra $\iak$ was represented as a
quotient of the universal differential envelope $\ude(\aaa)$ over
the simplicial differential ideal $\iii$ \rf{e58}. The form of the
K\"ahler differential for $\iak$ is given in \rf{e62}.

\section{Functorial properties}

As it was already mentioned, an arbitrary monotone mapping $\pkk$ 
between two posets $\kkk$ and $\kkk'$ produces no homomorphism of 
their incidence algebras. However, if we consider the category 
$\SC$ of simplicial complexes, the situation becomes completely 
different. In this section I show that that the correspondence 
$\kkk\mapsto\iak$ is a contravariant functor from the category 
$\SC$ to the category $\DMCA$ of differential moduli over 
commutative algebras.

\paragraph{The category $\DMCA$.} The objects of the category 
$\DMCA$ are differential moduli over semisimple commutative 
algebras --- see Definition \ref{dfdm}. The morphisms of $\DMCA$ 
are {\em differentiable mappings}.

\begin{dfn} Let $\dm=(\ia,\aaa,\kd)$, $\dm'=(\ia',\aaa',\kd')$ be
two differential moduli. A mapping $\phi:\ia\to\ia'$ is called
{\dff differentiable} iff

\be{edmca}
\begin{array}{l}
\phi\quad\text{\em is a homomorphism of graded algebras}\\
\phi\circ\kd'=\kd\circ\phi
\end{array}
\end{equation}

\end{dfn}

\begin{lemma} Any differentiable mapping is completely defined by 
its values on the basic algebra $\aaa$.

\end{lemma}

\begin{proof}
By induction, let $w\in\ia^1$, then $w=\sum_ia_i\,\kd{}b_i$ with
$a_i,b_i\in\aaa$. Then
\(\phi(w)=\sum\phi(a_i)\phi(\kd{}b_i)=\sum\phi(a_i)\kd'\phi(b_i)\).
When $w^{n+1}\in\ia^{n+1}$, it reads
\(w^{n+1}=\sum_ia_i\,\kd{}w^n_i\) with $a_i\in\aaa$ and
$w^n_i\in\ia^n$. Applying of $\phi$ and using the second property
\rf{edmca} completes the proof.
\end{proof}

So, the notion of differentiable mappings can be referred to
homomorphisms $\phi:\aaa\to\aaa'$ between basic algebras. Let
$\dm=(\ia,\aaa,\kd)$, $\dm'=(\ia',\aaa',\kd')$ be two differential
moduli. A mapping $\phi:\aaa\to\aaa'$ completely determines both
the homomorhisms $\phi:\ia\to\ia$ and $\phio:\ude\to\ude$. Then
$\phi$ is differentiable if and only if the following holds:

\be{edi}
\phio(\iii)\subseteq\iii'
\end{equation}

Let the basic algebras $\aaa,\aaa'$ are represented by functions on
the sets $\kkk$, $\kkk'$, respectively. Then the dual mapping $\pkk$
determines $\phi$ completely, and the explicit form of the mapping
$\phio$ is the following:  for any $\story{P_0\cdt P_n}\in\ude^n$

\be{ephio}
\phi\story{P_0\cdt P_n}\:=
\sum_{P_i':\,\forall i\,P_i'\pp=p_i}
\story{P_0'\cdt P_n'}
\end{equation}

\paragraph{Functorial properties.} To prove that the correspondence 
$\kkk\mapsto\iak$ is a functor we have to provide a correspondence 
between the arrows of the categories. As stated above, any 
(set-theoretical) mapping $\pkk$ gives rise to a homomorphism 
\rf{ephio} $\phi=\pp^*$ of the basic algebras.

\begin{theorem}
Let $\kkk$, $\kkk'$ be two simplicial complexes, and let a mapping
$\pkk$ be simplicial. Then its dual $\phi=\pp^*$ is differentiable.
\end{theorem}

\begin{proof}
By virtue of \rf{edi} it suffices to prove that simplicial ideals
are mapped into simplicial ideals in the target algebra. Let
$W'=\story{P_0'\cdt P_n'}$ be a fair story. Let us prove that its
image

\[
W'\pp=\story{P_0'\pp\cdt P_n'\pp}
\]

\noindent being a sequence of elements of $\kkk$ is either a fair
story or not a story. $P_i'=P_{i-1}'+v_i$ for any $i=1\cdt n$, and
$\forall i\neq j\;v_i\neq v_j$. Denote $P_i=P_i'\pp$ and
$v_i=v_i'\pp$ (recall that $\vvv'\pp\subseteq\vvv$). For each $i$
we have exactly two possibilities: either $v_i\in P_{i-1}$ or
$v_i\not\in P_{i-1}$ (and then according to \rf{e80}
$P_i=P_{i-1}+v_i\in\kkk$). If $v_i\in P_i$ for some for some $i$
then $W'\pp$ contains $P_i=P_{i-1}$ (and therefore is not a story)
otherwise $\forall i\; P_i=P_{i-1}+v_i$, so $W$ is a fair story.
So,

\be{efsm}
\{\text{\texttt{fair stories}}\}\pp\subseteq
\{\text{\texttt{fair stories}}\} 
\end{equation}

Return to simplicial ideals. Let $\story{P_0\cdt P_n}$ be an unfair
story, then $\phi(\story{P_0\cdt P_n})$ is according to \rf{ephio}
a sum of unfair stories (otherwise \rf{efsm} is violated).

Now let $i=\epsilon_w-\epsilon_{\tilde{w}}\in\iii_S$ such that
$P_n=P_0+v_1+\dotsb+v_n$ in $w$ and
$P_n=P_0+v_{j_1}+\dotsb+v_{j_n}$ in $\tilde{w}$. Let $w'=P_0'\cdt
P_n'$ be a fair story in $\kkk'$ such that $w'\pp=w$, then
necessarily

\begin{itemize}
\item \(P_n'=P_0'+v_1'+\dotsb+v_n'\)
\item all $v_j'$ are disjoint
\item \(\forall j\: v_j'\pp=v_j\) --- there is 1--1 correspondence
\end{itemize}

Make a new fair story $\tilde{w}'$ from $w'$ performing the same
permutation of vertices as that making $\tilde{w}$ from $w$, and we
obtain for it $\tilde{w}'\pp=\tilde{w}$ and
\(\epsilon_{\tilde{w}'}=\epsilon_{w'}\). So, any fair preimage of
$i\in\iii_S$ will be in $\iii_S'$

\end{proof}

\paragraph{Summary of this section.} The correspondence 
$\kkk\mapsto\iak$ is proved to be a contravariant functor from the 
category \(\SC=\{\text{\texttt{simplicial complexes}},\, 
\text{\texttt{simplicial mappings}}\}\) into the category $\DMCA$ 
of differential moduli over finite-dimensional semisimple 
commutative algebras.

\section{Concluding remarks}

It was shown that the incidence algebras of simplicial complexes
possess the structure of differential moduli, and it was
established that, contrary to posets of general form, the
correspondence between simplicial complexes and their incidence
algebras is a contravariant functor.

It occurs that simplicial complexes possess the natural structure
of discrete differential manifolds --- finite sets equiped with
differential calculi --- see, e.g. \cite{bdmh} for a review.

Aside of purely mathematical context, the presented results have
physical application being a basis for discrete approximations of
spacetime structure \cite{ijtp2k01}. When simplicial complexes are
treated as coarse-grained spacetime patterns, our results enable
the possibility of linking more and more refined approximations
between each other.

\paragraph{Ackowledgments.} The author appreciates numerous remarks
made by the participants of Friedmann Seminar on theoretical
physics (St.~Petersburg, Russia) headed by A.A.~Grib. Deep and
profound discussiond with my Italian colleagues G.~Landi, 
F.~Lizzi and the members of the Theory Group of the INFN Section in 
Naples were helpful.  Dr.  I.~Raptis (Pretoria) made his valuable 
impact explaining the physical meaning of the construction.

A financial support from the state research program `Universities
of Russia' is appreciated.

\end{document}